\documentclass[11pt]{amsart}

\usepackage[all]{xy}
\usepackage{amscd}
\usepackage{mathrsfs}
\usepackage{amssymb}
\usepackage{amsthm}
\usepackage{amsmath}

\newtheorem{theorem}{Theorem}[section]
\newtheorem{lemma}[theorem]{Lemma}
\newtheorem{proposition}[theorem]{Proposition}
\newtheorem{corollary}[theorem]{Corollary}
\newtheorem{remark}[theorem]{Remark}

\theoremstyle{definition}
\newtheorem{definition}[theorem]{Definition}
\newtheorem{example}[theorem]{Example}

\theoremstyle{remark}

\numberwithin{equation}{section}

\newcommand{\bde}{\begin{definition}}
\newcommand{\ede}{\end{definition}}
\newcommand{\bpr}{\begin{proposition}}
\newcommand{\epr}{\end{proposition}}
\newcommand{\bth}{\begin{theorem}}
\newcommand{\ethm}{\end{theorem}}
\newcommand{\bexa}{\begin{example}}
\newcommand{\eexa}{\end{example}}
\newcommand{\bcor}{\begin{corollary}}
\newcommand{\ecor}{\end{corollary}}
\newcommand{\blem}{\begin{lemma}}
\newcommand{\elem}{\end{lemma}}
\newcommand{\brem}{\begin{remark}}
\newcommand{\erem}{\end{remark}}
\newcommand{\bprf}{\begin{proof}}
\newcommand{\eprf}{\end{proof}}
\def \benu{\begin{enumerate}\renewcommand{\labelenumi}{(\roman{enumi})}\renewcommand{\itemsep}{0pt}}
\newcommand{\eenu}{\end{enumerate}}
\newcommand{\N}{\mathbb{N}}

\newcommand{\Z}{\mathbb{Z}}
\newcommand{\Q}{\mathbb{Q}}
\newcommand{\F}{\mathbb{F}}
\newcommand{\G}{\mathbb{G}}

\def\Coker{{\rm Coker\,}}

\def\Gal{{\rm Gal}}

\def\Pic{{\rm Pic}}

\def\Spec{{\rm Spec}}

\newcommand{\Hom}{{\rm Hom}}

\newcommand{\bsc}{\begin{tobira}}
\newcommand{\esc}{\end{tobira}}
\newcommand{\abc}{\xymatrix{0 \ar@{>}[r] & A \ar@{>}[r]^{f} & B \ar@{>}[r]^{g} & C \ar@{>}[r] & 0}}
\newcommand{\xyz}{\xymatrix{0 \ar@{>}[r] & X \ar@{>}[r]^{f} & Y \ar@{>}[r]^{g} & Z \ar@{>}[r] & 0}}
\def\phi{\varphi}
\def\hat{\widehat}
\def\lra{\longrightarrow}

\def\ra{\rightarrow}

\def\ol{\overline}

\def\ul{\underline}

\def\invlim#1{\lim\limits_{\substack{\longleftarrow\\#1}}}
\def\dirlim#1{\lim\limits_{\substack{\longrightarrow\\#1}}}
\DeclareMathOperator{\chara}{\mathrm{char}}

\begin{document}

\title[Curves over local fields]{Abelian \'etale coverings of modular curves over local fields}

\address{Harvard University, Department of Mathematics, 1 Oxford Street, Cambridge, MA 02138, USA}
\email{yoshida\char`\@math.harvard.edu}

\date{\today}

\begin{abstract}
We relate a part of the abelian \'etale fundamental group of curves over local fields to the component group of the N\'eron model of the jacobian. We apply the result to the modular curve $X_0(p)/\Q_p$ to show that the unramified abelian covering $X_1(p) \ra X_0(p)$ (Shimura covering) uses up all the possible ramification over the special fiber of $X_0(p)$. 
\end{abstract}

\maketitle

\tableofcontents

\section{Introduction}

The geometric class field theory tries to classify the abelian \'etale coverings of a proper smooth variety \cite{KatzLang}, and there is an arithmetic interest when the base field is a number field or a local field. There is a well-formulated class field theory for these curves over global/local fields (e.g.\ Kato-Saito \cite{KatoSaito}, Saito \cite{SSaito}) but there are very few examples where the abelian \'etale coverings are explicitly classified. The most remarkable known example is the case of modular curve $X_0(p)$ over $\Q$ (Mazur \cite{Mazur}, Introduction, Theorem (2)), which makes use of the deep theory of integral Hecke algebras. We give the corresponding result for $X_0(p)$ over $\Q_p$, by proving general results for the curves over local fields.

For a proper variety $X$ over a field $K$, let $\pi_1^{ab}(X)$ be the abelian \'etale fundamental group of $X$. There is a natural surjection $\pi_1^{ab}(X)\ra G_K^{ab}$ where $G_K^{ab}=\Gal(K_{ab}/K)$ is the Galois group of the maximal abelian extension of $K$, and we denote the kernel by $\pi_1^{ab}(X)^{\text{\rm geo}}$:
\[\xymatrix{0 \ar[r] & \pi_1^{ab}(X)^{\text{\rm geo}} \ar[r] & \pi_1^{ab}(X) \ar[r] & G_K^{ab} \ar[r] & 0}\]
For a proper smooth geometrically irreducible variety $X$ over a {\it local field} $K$, i.e.\ a complete discrete valuation field $K$ with finite residue field $F$ with $\chara F=p$, we showed in our previous paper  \cite{TY} that $\pi_1^{ab}(X)^{\text{\rm geo}}$ is an extension of $\hat{\Z}^r$ by a finite torsion group $\pi_1^{ab}(X)^{\text{\rm geo}}_{\text{\rm tor}}$. Here the rank $r$ is the $F$-rank (the dimension of the maximal $F$-split subtorus) of the special fiber of the N\'eron model of the Albanese variety of $X$.

In this paper, we confine ourselves to the case of {\it curves}, and investigate the most mysterious part of $\pi_1^{ab}(X)^{\text{\rm geo}}_{\text{\rm tor}}$, namely $\pi_1^{ab}(X)^{\text{\rm geo}}_{\text{\rm ram}}$ which classifies the abelian \'etale coverings of the generic fiber which are ``completely ramified over the special fiber", defined as follows. Let $X_F$ be the special fiber of the minimal regular model ${\mathscr X}$ over the integer ring $O_K$ (Abyhankar \cite{Ab}). Then we have a natural map $\pi_1^{ab}(X)^{\text{\rm geo}}\lra \pi_1^{ab}(X_F)^{\text{\rm geo}}$ which is surjective if $X$ has a $K$-rational point, and denote the kernel by $\pi_1^{ab}(X)^{\text{\rm geo}}_{\text{\rm ram}}$:
\[\xymatrix{0 \ar[r] & \pi_1^{ab}(X)^{\text{\rm geo}}_{\text{\rm ram}} \ar[r] & \pi_1^{ab}(X)^{\text{\rm geo}} \ar[r] & \pi_1^{ab}(X_F)^{\text{\rm geo}} \ar[r] & 0}\]
Then our main result is:

\begin{theorem} \label{mainthm}
Assume that $X$ admits a $K$-rational point.
\benu
\item {\rm (Theorem \ref{georamell})} The dual of the prime-to-$p$ part of the finite abelian group $\pi_1^{ab}(X)^{\text{\rm geo}}_{\text{\rm ram}}$ injects to the component group $\Phi$ of the N\'eron model of the jacobian variety of $X$.
\item {\rm (Theorem \ref{ppartprop})} Assume moreover that the absolute ramification index $e$ of $K$ is less than $p-1$, $\chara K=0$, and $X$ has semistable reduction over $O_K$. Then the $p$-primary part of $\pi_1^{ab}(X)^{\text{\rm geo}}_{\text{\rm ram}}$ vanishes.
\eenu
\end{theorem}

Note that if $X_F$ is smooth, i.e.\ $X$ has good reduction, $\pi_1^{ab}(X)^{\text{\rm geo}}\ra \pi_1^{ab}(X_F)^{\text{\rm geo}}$ is isomorphic on the prime-to-$p$ part by the proper smooth base change theorem on $H^1_{et}$. As $\pi_1^{ab}(X_F)^{\text{\rm geo}}$ can be calculated in principle from the special fiber $X_F$, we have control on the whole of $\pi_1^{ab}(X)$ in the case $e<p-1$, $\chara K=0$ and $X$ has semistable reduction, in terms of the special fiber of the N\'eron model of the jacobian of $X$. 

In the latter part of this paper, we give an application of these results to the modular curve $X_0(p)$, where $X_0(p)$ is the usual modular curve classifying the elliptic curves with $\Gamma_0(p)$-structures ($p$ is a prime). For the curve $X_0(p)/\Q_p$, above theorem enables us to compute the group $\pi_1^{ab}(X_0(p)/\Q_p)^{\text{\rm geo}}$ completely:

\bth \label{X0pcompleteintro}
$\pi_1^{ab}(X_0(p)/\Q_p)^{\text{\rm geo}}$ has the following structure:
\[\xymatrix{0 \ar[r] & \Phi(J_0(p)) \ar[r] & \pi_1^{ab}(X_0(p)/\Q_p)^{\text{\rm geo}} \ar[r] & \hat{\Z}^r \ar[r] & 0}\]
where $\Phi(J_0(p))\cong \pi_1^{ab}(X_0(p)/\Q_p)^{\text{\rm geo}}_{\text{\rm ram}}=\pi_1^{ab}(X_0(p)/\Q_p)^{\text{\rm geo}}_{\text{\rm tor}}$ is a cyclic group of order equal to the numerator of $\frac{p-1}{12}$, and $r=\frac{g+h-1}{2}$ where $g$ is the genus of $X_0(p)$ and $h$ is the number of the supersingular points defined over $\F_p$.
\ethm

{\bf Acknowledgements.} This paper was mostly written in 2001, as a part of the author's master thesis at the University of Tokyo. The author would like to express his sincere gratitude to his thesis adviser K. Kato for suggesting the problem and his constant encouragement. He would like to thank K.\ Ban for helpful discussions, and B.\ Conrad for many comments on the first draft.

{\bf Notations.} Throughout this paper, a {\it local field} $K$ means a complete discrete valuation field with finite residue field $F$, with $\chara F=p$. $O_K$ is the integer ring of $K$. For any field $K$, $\ol{K}$ is a separable closure of $K$, and $G_K=\Gal(\ol{K}/K)$ is the absolute Galois group of $K$. For a variety $X$ over any field $K$, $\ol{X}=X\times_K \Spec(\ol{K})$. 

For an abelian group $X$, its Pontrjagin dual is denoted by $X^\vee=\Hom(X,\Q/\Z)$. $\hat{\Z}\cong \prod\limits_p\Z_p$ is the profinite completion of $\Z$. 

For a scheme $X$ over $O_K$, we denote the generic fiber and the special fiber of $X$ respectively by $X_K,X_F$. For any group scheme $X$, $X[m]$ denotes the kernel of the multiplication-by-$m$ map, which is also a group scheme. We use the notation $\mu_N=\G_m[N]$ for the group scheme of $N$-th roots of unity over an arbitrary base, without specifying the base scheme.

\section{Preliminaries on abelian \'etale fundamental groups}

Here we review some generalities on abelian \'etale fundamental groups, and fix the notations. 

\subsection{Geometric abelian \'etale fundamental groups}

For any noetherian connected scheme $X$, the \'etale fundamental group $\pi_1(X)$ is a profinite group classifying finite \'etale coverings of $X$ (\cite{SGA1}), and we denote the maximal abelian quotient by $\pi_1^{ab}(X)$. When $X$ is a proper variety over a field $K$, there is an exact sequence ({\it ibid.}, Expos\'e IX, Th.\ 6.1):
\begin{equation} \label{pi1}
\xymatrix{1 \ar[r] & \pi_1(\ol{X}) \ar[r] & \pi_1(X) \ar[r] & G_K \ar[r] & 1}
\end{equation}
where $\ol{K}$ is a separable closure of $K$, $\ol{X}=X\times_K\Spec(\ol{K})$, and $G_K=\Gal(\ol{K}/K)$ is the absolute Galois group of $K$. Therefore we have surjection $\pi_1^{ab}(X)\ra G_K^{ab}$ where $G_K^{ab}=\Gal(K_{ab}/K)$ is the Galois group of the maximal abelian extension of $K$, and we denote the kernel by $\pi_1^{ab}(X)^{\text{\rm geo}}$:
\begin{equation} \label{pi1ab}
\xymatrix{0 \ar[r] & \pi_1^{ab}(X)^{\text{\rm geo}} \ar[r] & \pi_1^{ab}(X) \ar[r] & G_K^{ab} \ar[r] & 0}
\end{equation}

\begin{remark} \label{rat}
{\rm When $X$ has a $K$-rational point $x$, (\ref{pi1}) and consequently also (\ref{pi1ab}) has a splitting, and $\pi_1^{ab}(\ol{X})_{G_K}\cong \pi_1^{ab}(X)^{\text{\rm geo}}$ where the former group is the coinvariant with respect to the $G_K$-action by inner automorphisms. In this case, $\pi_1^{ab}(X)^{\text{\rm geo}}$ has a geometric interpretation as the group classifying the abelian finite \'etale coverings of $X$ in which $x$ splits completely (\cite{KatzLang}). }
\end{remark}

Now assume that $K$ is a {\it local field}, i.e.\ a complete discrete valuation fields with a finite residue field $F$ with $\chara F=p$, and let ${\mathscr X}$ be a proper flat model of $X$ over the integer ring $O_K$, i.e. a proper flat scheme ${\mathscr X}\ra \Spec\ O_K$ with ${\mathscr X}\times_{O_K}\Spec\ K\cong X$, and let $X_F={\mathscr X}\times_{O_K}\Spec\ F$ be the special fiber. Then we have a canonical surjection:
\[\pi_1^{ab}(X)\lra \pi_1^{ab}({\mathscr X})\cong \pi_1^{ab}(X_F)\]
Here the latter isomorphism follows by  \cite{SGA1} Expos\'e X, Th.\ 2.1 (a part of proper base change theorem for $H^1_{et}$). We have the corresponding map on the geometric part by the commutative diagram:
\begin{equation} \label{cdpi1}
\xymatrix{0 \ar[r] & \pi_1^{ab}(X)^{\text{\rm geo}} \ar[d]\ar[r] & \pi_1^{ab}(X) \ar[d]\ar[r] & G_K^{ab} \ar[d]\ar[r] & 0\\
0 \ar[r] & \pi_1^{ab}(X_F)^{\text{\rm geo}} \ar[r] & \pi_1^{ab}(X_F) \ar[r] & G_F \ar[r] & 0}
\end{equation}

\bde
Denote the kernel of left and middle vertical maps of (\ref{cdpi1}) respectively by $\pi_1^{ab}(X)^{\text{\rm geo}}_{{\mathscr X}\text{-ram}}$, $\pi_1^{ab}(X)_{{\mathscr X}\text{-ram}}$. (These groups classify the abelian \'etale coverings of $X$ which become completely ramified over a part of $X_F$ when extended to the coverings over $X_F$.)
\ede

By definition, we have an exact sequence:
\begin{equation} \label{georam}
\xymatrix{0 \ar[r] & \pi_1^{ab}(X)^{\text{\rm geo}}_{{\mathscr X}\text{-ram}} \ar[r] & \pi_1^{ab}(X)^{\text{\rm geo}} \ar[r] & \pi_1^{ab}(X_F)^{\text{\rm geo}}}
\end{equation}

As $\pi_1^{ab}(\ol{X})$ is the Pontrjagin dual of the \'etale cohomology group $H^1_{et}(\ol{X},\Q/\Z)$, we can pass to the dual and translate the above definition into the language of $H^1_{et}$. By the Hochschild-Serre spectral sequece, we have a short exact sequence:
\[\xymatrix{0 \ar[r] & H^1(G_K,\Q/\Z) \ar[r] & H^1_{et}(X,\Q/\Z) \ar[r] & H^1_{et}(\ol{X},\Q/\Z)^{G_K} \ar[r] & 0}\]
where $(-)^{G_K}$ is the Galois invariant. (Here we used $H^2(G_K,\Q/\Z)=0$ for local fields $K$, and the same holds when we replace $K$ by the finite field $F$.) We can view this exact sequence as the Pontrjagin dual of (\ref{pi1ab}), and therefore we have a canonical isomorphism:
\begin{equation} \label{pi1geoasH1}
(\pi_1^{ab}(X)^{\text{\rm geo}})^\vee\cong H^1_{et}(\ol{X},\Q/\Z)^{G_K}
\end{equation}
Hence by taking the Pontrjagin dual of the exact sequence (\ref{georam}), we have:

\blem \label{firstlem}
Let $X$ be a proper geometrically irreducible variety over a local field $K$, and let ${\mathscr X}$ be its proper flat model over $O_K$, and denote the special fiber by $X_F$. Then there is a canonical exact sequence of abelian groups:
\begin{equation} \label{firstlemseq}
\xymatrix{H^1_{et}(\ol{X_F},\Q/\Z)^{G_F} \ar[r] & H^1_{et}(\ol{X},\Q/\Z)^{G_K} \ar[r] & (\pi_1^{ab}(X)^{\text{\rm geo}}_{{\mathscr X}\text{-ram}})^\vee \ar[r] & 0}
\end{equation}
\elem

\brem \label{ratpt}
{\rm Note that when $X$ has a $K$-rational point, it gives a splitting of the exact sequence (\ref{pi1ab}), and consequently all the vertical arrows of (\ref{cdpi1}) are surjective, and exact sequences (\ref{georam}) and (\ref{firstlemseq}) turn out to be the short exact sequences.}
\erem

\subsection{\'Etale cohomology groups}

Here we review some of the properties of the group $H^1_{et}(\ol{X},\Q/\Z)^{G_K}$ for a proper geometrically irreducible {\it curve} $X$ over a field $K$. 

First recall that, for prime-to-$p$ part, we have:

\begin{lemma}[\cite{SGA4}, Expos\'e IX, Corollaire 4.7] \label{Picell}
Let $\ell$ be a prime different from $\chara K$. For a proper geometrically connected curve $X$ over a field $K$, $H^1_{et}(\ol{X},\Q_\ell/\Z_\ell)\cong \Pic^0(X)[\ell^\infty](-1)$ as $G_K$-modules.
\end{lemma}

More generally, considering the long exact sequence induced from the short exact sequence:
\[\xymatrix{
0 \ar[r] & \frac1N\Z/\Z \ar[r] & \Q/\Z \ar[r]^-N & \Q/\Z \ar[r] & 0
}\]
we have $H^1_{et}(\ol{X},\frac1N\Z/\Z)=H^1_{et}(\ol{X},\Q/\Z)[N]$. As $\Q/\Z=\dirlim{N}\frac1N\Z/\Z$ and the \'etale cohomology commutes with the directed inductive limits (\cite{SGA4} Expos\'e VII, Th\'eor\`eme 5.7), we deduce that:
\[H^1_{et}(\ol{X},\Q/\Z)=\dirlim{N}H^1_{et}(\ol{X},\frac1N\Z/\Z)=\bigcup_NH^1_{et}(\ol{X},\frac1N\Z/\Z)\]
(In particular, $H^1_{et}(\ol{X},\Q/\Z)$ is an torsion module.) Therefore we have:
\begin{equation} \label{Q/ZtoN}
H^1_{et}(\ol{X},\Q/\Z)^{G_K} = \bigcup_NH^1_{et}(\ol{X},\Q/\Z)[N]^{G_K}=\bigcup_NH^1_{et}(\ol{X},\frac1N\Z/\Z)^{G_K}
\end{equation}
and can reduce the calculation of $H^1_{et}(\ol{X},\Q/\Z)^{G_K}$ to that of $H^1_{et}(\ol{X},\Z/N\Z)^{G_K}$, for which we can make use of the following generalization of Lemma \ref{Picell}:

\begin{lemma}[Milne \cite{Milne}, Theorem 3.9, Proposition 4.16] \label{PicN}
For any proper geometrically irreducible curve over a field $K$ and integer $N\geq 1$, we have a canonical isomorphism of $G_K$-modules:
\[H^1_{et}(\ol{X},\Z/N\Z)\cong \Hom_{\rm gp}(\mu_N,\Pic^0(X)[N])\]
where $\Hom_{\rm gp}$ denotes the $G_K$-module consisting of homomorphisms as group schemes. In particular, we have:
\[H^1_{et}(\ol{X},\Z/N\Z)^{G_K}\cong \Hom_{K\text{\rm -gp}}(\mu_N,\Pic^0(X)[N])\]
where $\Hom_{K\text{\rm -gp}}$ denotes the abelian group consisting of homomorphisms defined over $K$.
\elem

\brem
{\rm This interpretation shows that, $(\pi_1^{ab}(X)^{\text{\rm geo}})^\vee\cong H^1_{et}(\ol{X},\Q/\Z)^{G_K}$ is the {\it maximal $\mu$-type subgroup} of $\Pic^0(X)$ in the terminology of Mazur \cite{Mazur}, I-3. This re-interpretes the Mazur's theory of the $\mu$-type subgroup of the jacobian of the modular curves as the theory of $\pi_1^{ab}(X)^{\text{\rm geo}}$ of the modular curves.}
\erem

\section{The prime-to-$p$ part of $\pi_1^{ab}(X)_{\text{\rm ram}}^{\text{\rm geo}}$}

From this section, we assume that $X$ is {\it a proper smooth geometrically irreducible curve over a local field $K$}. In this case, the proper flat regular model ${\mathscr X}$ always exists by Abyhankar \cite{Ab}, and denote its special fiber by $X_F$. 

Let $J=\Pic^0(X)$ be the jacobian variety of $X$, and denote the N\'eron model and its special fiber respectively by ${\mathscr J}$ and $J_F$. The quotient of $J_F$ by its connected component of the identity $J_F^0$ is a finite \'etale group $\Phi(J)$ over $F$, which we call the {\it group of components} of ${\mathscr J}$ (or of $J$, by an abuse of language). 

Now we assume that $X$ {\it admits a $K$-rational point} (see Remark \ref{ratpt}), and our starting point of the investigation is the following:

\begin{lemma}[Raynaud \cite{Raynaud},  \cite{SGA7} Expos\'e IX, (12.1.12)] \label{Raylem}
Assume that $X$ admits a $K$-rational point. Then we have $\Pic^0(X_F)\cong J_F^0$ as smooth group schemes over $F$.
\elem

This lemma, in view of Lemma \ref{firstlem} and \ref{PicN}, assures that $\pi_1^{ab}(X)^{\text{\rm geo}}_{{\mathscr X}\text{-ram}}$ is independent of the model ${\mathscr X}$. As we assume the existence of $K$-rational point in the rest of this paper, we will suppress the notation ${\mathscr X}$ from $\pi_1^{ab}(X)^{\text{\rm geo}}_{{\mathscr X}\text{-ram}}$. 


In this section we will treat the prime-to-$p$ part of $\pi_1^{ab}(X)_{\text{\rm ram}}^{\text{\rm geo}}$, and the main result is stated as follows:

\bth \label{georamell}
Assume that $X$ admits a $K$-rational point. Then there is an injection of finite abelian groups $(\pi_1^{ab}(X)^{\text{\rm geo}}_{\text{\rm ram}})^\vee_{{\rm not}\, p} \lra \Phi(J)$, where $(-)_{{\rm not}\, p}$ denotes the prime-to-$p$ part ($p=\chara F$).
\ethm

We will start by fixing some notations on Galois modules, and give the proof of the above theorem. As the theorem can be proved by establishing the injectivity on $\ell$-primary parts for each prime number $\ell\neq p=\chara F$, we fix an $\ell\neq p$ in the remainder of this section. 

\subsection{Galois modules}

For any smooth group scheme $X$ of finite type over $K$, define pro-$\ell$ (resp. ind-$\ell$) group scheme $T_\ell(X)$ (resp. $X[\ell^\infty]$) by:
\[T_\ell(X)=\invlim{n} X[\ell^n],\ \ X[\ell^\infty]=\dirlim{n} X[\ell^n]\]
where $X[m]$ denotes the kernel of the multiplication-by-$m$ map which is an \'etale group scheme over $K$, and inductive (resp.\ projective) limit is taken with respect to the canonical inclusions (resp.\ multiplication-by-$\ell$ maps). We use the same notations for group schemes over $O_K$ or $F$. We often identify $T_\ell(X),X[\ell^\infty]$ with the associated Galois modules, which are $\Z_\ell$-modules and satisfies:
\[X[\ell^\infty]\cong T_\ell(X)\otimes_{\Z_\ell}(\Q_\ell/\Z_\ell),\ \ \ \ T_\ell(X)=\Hom_{\Z_\ell}(\Q_\ell/\Z_\ell, X[\ell^\infty])\]

Now we fix the notations concerning Galois modules. Let $G=G_K\ \text{or}\ G_F$ be the absolute Galois group, and let $M$ be an arbitrary $\Z_\ell$-module with a continuous action of $G$. As usual, we write $\Z_\ell(1)=T_\ell(\G_m)$, and define the {\it Tate twist} $M(r)$ of $M$ for $\forall r\in \Z$ by:
\[M(r)=\begin{cases}
M\otimes_{\Z_\ell} \Z_\ell(r) & (r>0)\\
M & (r=0)\\
\Hom_{\Z_\ell}(\Z_\ell(-r),M) & (r<0)
\end{cases}
\]
as Galois modules. $M\mapsto M(r)$ gives exact functors of Galois modules for $\forall r\in \Z$, and we have canonical isomorphisms $M(r)(s)\cong M(r+s)$ for $\forall r,s\in \Z$. The Galois action on $\Z_\ell(1)$ is via the {\it cyclotomic character} $\chi:G\ra \Z_\ell^\times$. The Galois invariant and coinvariant is denoted respectively by $M^G,M_G$. 

\bde
The {\it $\chi$-part}\ $M^\chi$ of $M$ is defined by $M^\chi=(M(-1)^G)(1)$. It is canonically identified with the maximal subgroup of $M$ on which $G$ acts via the cyclotomic character $\chi$.
\ede

Note that the functors $M\mapsto M^G,\ M\mapsto M^\chi$ are both left exact. 

\subsection{Tate modules of the jacobians of curves}

Now we return to the curve $X$ over a local field $K$. Applying Lemma \ref{Picell} and Lemma \ref{Raylem}, we can restate the Lemma \ref{firstlem} in the present case, in view of Remark \ref{ratpt}, as follows:

\blem \label{firstlem2}
For a proper smooth geometrically irreducible curve $X$ over a local field $K$ which admits a $K$-rational point, there is a canonical short exact sequence of abelian groups:
\[\xymatrix{0 \ar[r] & J_F^0[\ell^\infty](-1)^{G_F} \ar[r] & J[\ell^\infty](-1)^{G_K} \ar[r] & (\pi_1^{ab}(X)^{\text{\rm geo}}_{\text{\rm ram}})^\vee_\ell \ar[r] & 0}
\]
where $(\pi_1^{ab}(X)^{\text{\rm geo}}_{\text{\rm ram}})^\vee_\ell$ denotes the $\ell$-primary part of $(\pi_1^{ab}(X)^{\text{\rm geo}}_{\text{\rm ram}})^\vee$.
\elem

Using the canonical isomorphism $M(-1)^G\cong M^\chi(-1)$ of $G_F$-modules, we can express this lemma by the short exact sequence of $G_F$-modules:
\begin{equation} \label{georamellses}
\xymatrix{0 \ar[r] & J_F^0[\ell^\infty]^\chi \ar[r] & J[\ell^\infty]^\chi \ar[r] & (\pi_1^{ab}(X)^{\text{\rm geo}}_{\text{\rm ram}})^\vee_\ell(1) \ar[r] & 0}
\end{equation}

Now we relate this group with the component group $\Phi(J)$. For this we need the description of $\Phi(J)$ by the Galois module associated to $J$, stated as follows:

\begin{lemma}[ \cite{SGA7} Expos\'e IX, Proposition 11.2]
For an abelian variety $A$ over a local field $K$, denote the N\'eron model and its special fiber respectively by ${\mathscr A},A_F$. Let $A_F^0$ be the connected component of the identity, and denote the component group by $\Phi(A)=A_F/A_F^0$. For any prime $\ell\neq \chara F$, there is a canonical isomorphism of $G_F$-modules:
\[\Phi(A)_\ell \cong \Coker (T_\ell(A)^I\otimes (\Q_\ell/\Z_\ell)\lra (T_\ell(A)\otimes (\Q_\ell/\Z_\ell))^I)\]\
where $\Phi(A)_\ell$ is the $\ell$-primary part of $\Phi(A)$, and $(-)^I$ denotes the invariant by the action of the inertia group $I\subset G_K$.
\elem

By the canonical isomorphism of $G_F$-modules $T_\ell(A)^I\cong T_\ell(A_F^0)$ (\cite{SGA7} Expos\'e IX, Proposition 2.2.5, (2.2.3.3)), we can restate the above lemma by the short exact sequence of $G_F$-modules:
\begin{equation} \label{phiell}
\xymatrix{0 \ar[r] & A_F^0[\ell^\infty] \ar[r] & A[\ell^\infty]^I \ar[r] & \Phi(A)_\ell \ar[r] & 0}
\end{equation}

\begin{proof}[Proof of Theorem \ref{georamell}]
Applying the left exact functor $M\mapsto M^\chi$ to the short exact sequence (\ref{phiell}) for $A=J$, we have the exact sequence:
\begin{equation}
\xymatrix{0 \ar[r] & J_F^0[\ell^\infty]^\chi \ar[r] & J[\ell^\infty]^\chi \ar[r] & \Phi(J)_\ell^\chi}
\end{equation}
Comparing this with the short exact sequence (\ref{georamellses}), we have a canonical injection $(\pi_1^{ab}(X)^{\text{\rm geo}}_{\text{\rm ram}})^\vee_\ell\lra \Phi(J)_\ell^\chi(-1)$. As $\Phi(J)_\ell^\chi(-1)\cong \Phi(J)_\ell(-1)^G$ injects to $\Phi(J)_\ell$ as finite abelian groups and $\ell\neq p$ was arbitrary, we have proven the theorem.
\eprf

\subsection{Example: elliptic curves}
Let $X=E$ be an elliptic curve over the local field $K$ (see Silverman \cite{Silverman}, Chap.\ IV). (\ref{georamellses}),(\ref{phiell}) for $A=E$ reads:
\begin{gather}
\xymatrix{0 \ar[r] & E_F^0[\ell^\infty]^\chi \ar[r] & E[\ell^\infty]^\chi \ar[r] & (\pi_1^{ab}(E)^{\text{\rm geo}}_{\text{\rm ram}})^\vee_\ell(1) \ar[r] & 0}\\
\xymatrix{0 \ar[r] & E_F^0[\ell^\infty] \ar[r] & E[\ell^\infty]^I \ar[r] & \Phi(E)_\ell \ar[r] & 0}
\end{gather}
When $E$ has good reduction over $O_K$, trivially $(\pi_1^{ab}(X)^{\text{\rm geo}}_{\text{\rm ram}})_{{\rm not}\, p}=0$. If $E$ has additive reduction, $E_F^0$ is a unipotent group, therefore $E_F^0[\ell^\infty]=0$. Hence we have canonical isomorphisms :
\[E[\ell^\infty]^\chi\cong (\pi_1^{ab}(E)^{\text{\rm geo}}_{\text{\rm ram}})^\vee_\ell(1),\ \ \ \ E[\ell^\infty]^I\cong \Phi(E)_\ell\]
which gives $\Phi(E)_\ell^\chi\cong (\pi_1^{ab}(E)^{\text{\rm geo}}_{\text{\rm ram}})^\vee_\ell(1)$. Moreover, in this case we know that the $G=G_F$ acts trivially on $\Phi(E)$, hence we have :

\bpr
If $E$ has good or additive reduction over $O_K$, the injection in Th.\ \ref{georamell} is a bijection $(\pi_1^{ab}(E)^{\text{\rm geo}}_{\text{\rm ram}})^\vee_{{\rm not}\, p} \cong \Phi(E)(-1)^G\cong \Phi(E)[q-1]$, where $q=|F|$.
\epr

In the multiplicative reduction case, this is not necessarily true.

\section{The $p$-primary part for low absolute ramification case}

As in the previous section, we consider a proper smooth geometrically irreducible curve $X$ over a local field $K$ which admits a $K$-rational point. Let $e\in \N$ denote the {\it absolute ramification index} of $K$, i.e.\ the normalized valuation of $p=\chara F$ in $K$. In this section, we prove the following:

\bth \label{ppartprop}
Assume $e<p-1$ and $\chara K=0$. If $X$ has semistable reduction over $O_K$, then $(\pi_1^{ab}(X)^{\text{\rm geo}}_{\text{\rm ram}})_p=0$, where $(-)_p$ denotes the $p$-primary part ($p=\chara F$).
\ethm

\brem
{\rm This theorem is a generalization of the Prop.\ 7 of Kato-Saito \cite{KatoSaito}, where it is proved in the good reduction case. Here we employ a completely different method of proof.}
\erem

Our task here is to analyze the $H^1_{et}(\ol{X},\Z/p^n\Z)^{G_K}\cong \Hom_{K\text{\rm -gp}}(\mu_{p^n},\Pic^0(X)[p^n])$ (Lemma \ref{PicN}). For this purpose, we recall the structure of the $p^n$-torsion points $J[p^n]$ of the jacobian $J=\Pic^0(X)$, following  \cite{SGA7}, Expos\'e IX. In the rest of this section, we assume that {\it $X$ has semistable reduction over $O_K$}. In this case, $J$ has semistable reduction, i.e. the connected component of the identity $J_F^0$ of the special fiber $J_F$ of the N\'eron model ${\mathscr J}$ is a semiabelian scheme over $F$. Moreover, if we denote by ${\mathscr J}^0$ the connected component of the identity of ${\mathscr J}$, we know that ${\mathscr J}^0[p^n]$ is a quasi-finite group scheme over $O_K$ (\cite{SGA7} Expos\'e IX, Lemme 2.2.1)

Now recall that any quasi-finite group scheme ${\mathscr G}$ over $O_K$ has the canonical and functorial decomposition ${\mathscr G}={\mathscr G}^f\coprod {\mathscr G}'$ where ${\mathscr G}^f$ ({\it the fixed part} of ${\mathscr G}$) is finite flat over $O_K$ and ${\mathscr G}'$ has empty special fiber (\cite{EGA} II,(6.2.6),  \cite{SGA7} Expos\'e IX, (2.2.3.1)). 

Coming back to our case, ${\mathscr J}^0[p^n]$ has the decomposition:
\[{\mathscr J}^0[p^n]={\mathscr J}^0[p^n]^f\coprod {\mathscr J}^0[p^n]'\]
where ${\mathscr J}^0[p^n]^f$ is a finite flat group scheme. We denote the generic fiber of ${\mathscr J}^0[p^n]^f$ by $J^0[p^n]^f$, which is a subgroup of $J[p^n]$, the generic fiber of ${\mathscr J}^0[p^n]$ (recall $({\mathscr J}^0)_K=J$). Note that ${\mathscr J}^0[p^n],\ {\mathscr J}^0[p^n]^f$ have a common special fiber $J_F^0[p^n]$, which is canonically isomorphic to $\Pic^0(X_F)[p^n]$ by Lemma \ref{Raylem} (recall the running hypothesis that $X$ has a $K$-rational point).

Moreover, we have additional information by the semi-stability hypothesis (\cite{SGA7} Expos\'e IX, (5.5.8)):
\begin{equation} \label{J/J0fetale}
J[p^n]/J^0[p^n]^f\cong M_K\otimes (\Z/p^n\Z)=(M\otimes (\Z/p^n\Z))_K
\end{equation}
where $M$ denotes the character group of the toric part of $J_F$, which is an unramified Galois module, free of finite rank over $\Z$. We consider $M$ as an \'etale group scheme over $O_K$, therefore justifying the notation $M_K$. (Note that in general we have to take the corresponding object for the dual abelian variety, but we have the autoduality of the jacobian in the present case.)

Now we apply the result of Raynaud \cite{Raynaudp} to obtain the following lemma: 

\blem \label{ptorsisom}
Under the hypothesis $e<p-1$ and $\chara K=0$, the following natural homomorphism is an isomorphism:
\[\xymatrix{H^1_{et}(\ol{X_F},\Z/p^n\Z)^{G_F} \ar[r]^-\cong & H^1_{et}(\ol{X},\Z/p^n\Z)^{G_K}}\]
\elem

\bprf
By Lemma \ref{PicN}, we have:
\begin{align*}
H^1_{et}(\ol{X_F},\Z/p^n\Z)^{G_F} &\cong \Hom_{F\text{\rm -gp}}(\mu_{p^n},J_F^0[p^n])\\
H^1_{et}(\ol{X},\Z/p^n\Z)^{G_K} &\cong \Hom_{K\text{\rm -gp}}(\mu_{p^n},J[p^n])
\end{align*}
By (\ref{J/J0fetale}), we have an exact sequence:
\begin{gather*}\xymatrix{
0 \ar[r] & \Hom_{K\text{\rm -gp}}(\mu_{p^n},J^0[p^n]^f) \ar[r] & \Hom_{K\text{\rm -gp}}(\mu_{p^n},J[p^n])}\\
\xymatrix{\ar[r] & \Hom_{K\text{\rm -gp}}(\mu_{p^n},(M\otimes (\Z/p^n\Z))_K)
}\end{gather*}
But by Corollaire 3.3.6 of Raynaud \cite{Raynaudp}, we have:
\[\Hom_{K\text{\rm -gp}}(\mu_{p^n},(M\otimes (\Z/p^n\Z))_K)\cong \Hom_{O_K\text{\rm -gp}}(\mu_{p^n},M\otimes (\Z/p^n\Z))=0\]
which vanishes simply because $\mu_{p^n}$ is connected and $M\otimes (\Z/p^n\Z)$ is \'etale. Therefore we have, using  \cite{Raynaudp}, Corollaire 3.3.6 again,
\[\Hom_{K\text{\rm -gp}}(\mu_{p^n},J[p^n]) \cong \Hom_{K\text{\rm -gp}}(\mu_{p^n},J^0[p^n]^f) \cong \Hom_{O_K\text{\rm -gp}}(\mu_{p^n},{\mathscr J}^0[p^n]^f)\]
Therefore it suffices to see that $\Hom_{O_K\text{\rm -gp}}(\mu_{p^n},{\mathscr J}^0[p^n]^f)\cong \Hom_{F\text{\rm -gp}}(\mu_{p^n},J_F^0[p^n])$, which is equivalent by Cartier duality (\cite{SGA3} Expos\'e VIIA, (3.3.1)) to:
\[\Hom_{O_K\text{\rm -gp}}({\mathscr J}^0[p^n]^{f*},\Z/p^n\Z)\cong \Hom_{F\text{\rm -gp}}(J_F^0[p^n]^*,\Z/p^n\Z)\]
where $^*$ denotes the Cartier dual ${\mathscr G}^*={\rm \ul{Hom}}({\mathscr G},{\mathbb G}_m)$. But as $\Z/p^n\Z$ is \'etale, this is clear by the Hensel's lemma (e.g.  \cite{EGA} IV, (18.5.12)).
\eprf

\begin{proof}[Proof of Theorem \ref{ppartprop}]
By Lemma \ref{firstlem}, it is enough to show that the canonical homomorphism $H^1_{et}(\ol{X_F},\Q_p/\Z_p)^{G_F}\lra H^1_{et}(\ol{X},\Q_p/\Z_p)^{G_K}$ is an isomorphism. But by exactly the same argument as in (\ref{Q/ZtoN}) we have:
\begin{gather*}
H^1_{et}(\ol{X_F},\Q_p/\Z_p)^{G_F}=\bigcup_nH^1_{et}(\ol{X_F},\Q_p/\Z_p)[p^n]^{G_F}=\bigcup_nH^1_{et}(\ol{X_F},\frac{1}{p^n}\Z/\Z)^{G_F}\\
H^1_{et}(\ol{X},\Q_p/\Z_p)^{G_K}=\bigcup_nH^1_{et}(\ol{X},\Q_p/\Z_p)[p^n]^{G_K}=\bigcup_nH^1_{et}(\ol{X},\frac{1}{p^n}\Z/\Z)^{G_K}
\end{gather*}
Therefore the desired isomorphism is deduced from the corresponding result in the each $p$-power torsion level, namely the Lemma \ref{ptorsisom}.
\eprf

Combining with Th.\ \ref{georamell}, we have:

\bcor \label{loweinj}
Assume $e<p-1$ and $\chara K=0$. If $X$ admits a $K$-rational point and has semistable reduction over $O_K$, then there is an injection of finite abelian groups $(\pi_1^{ab}(X)^{\text{\rm geo}}_{\text{\rm ram}})^\vee \lra \Phi(J)$.
\ecor

\brem
{\rm We have the canonical perfect dualiy:
\[\Phi\otimes_\Z \Phi\lra \Q/\Z\]
deduced by the autoduality of the jacobian and  \cite{SGA7} Expos\'e IX, (11.4.1) or Conjecture 1.3 (proven in this case). This gives the canonical identification $\Phi\cong \Phi^\vee$, which enables us to state Cor.\ \ref{loweinj} as the {\it surjection} $\Phi(J)\lra \pi_1^{ab}(X)^{\text{\rm geo}}_{\text{\rm ram}}$.}
\erem

\section{The case of the modular curves}

In this section, we apply the results of preceding section to the modular curve $X_0(p)/\Q_p$ for a prime $p$, which is a proper smooth geometrically irreducible curve over $\Q_p$. For the definition and basic properties of the modular curve $X_0(p)$, we refer to Mazur  \cite{Mazur}. The result of this section is summarized as follows:

\bth \label{X0pcomplete}
$\pi_1^{ab}(X_0(p)/\Q_p)^{\text{\rm geo}}$ has the following structure:
\[\xymatrix{0 \ar[r] & \Phi(J_0(p)) \ar[r] & \pi_1^{ab}(X_0(p)/\Q_p)^{\text{\rm geo}} \ar[r] & \hat{\Z}^r \ar[r] & 0}\]
where $\Phi(J_0(p))\cong \pi_1^{ab}(X_0(p)/\Q_p)^{\text{\rm geo}}_{\text{\rm ram}}=\pi_1^{ab}(X_0(p)/\Q_p)^{\text{\rm geo}}_{\text{\rm tor}}$ is a cyclic group of order equal to the numerator of $\frac{p-1}{12}$, and $r=\frac{g+h-1}{2}$ where $g$ is the genus of $X_0(p)$ and $h$ is the number of the supersingular points defined over $\F_p$.
\ethm

\brem
{\rm This result should be interpreted as the local analogue of the theorem of Mazur which asserts that $\pi_1^{ab}(X_0(p)/\Q)^{\text{\rm geo}}\cong \Phi(J_0(p))$ (\cite{Mazur}, Introduction, Theorem (2)). The fact that $\pi_1^{ab}(X_0(p)/\Q)^{\text{\rm geo}}$ is equal to $\pi_1^{ab}(X_0(p)/\Q_p)^{\text{\rm geo}}_{\text{\rm ram}}=(\pi_1^{ab}(X_0(p)/\Q_p)^{\text{\rm geo}}_{\text{\rm tor}}$ shows that the maximal abelian \'etale covering of $X_0(p)/\Q$ ``uses up" all the ramification allowed at the special fiber at $p$. The origin of this phenomenon remains to be clarified.}
\erem

\subsection{The part $\pi_1^{ab}(X_0(p)/\Q_p)^{\text{\rm geo}}_{\text{\rm ram}}$}

First, we apply the results of the preceding sections to determine the part $\pi_1^{ab}(X_0(p)/\Q_p)^{\text{\rm geo}}_{\text{\rm ram}}$. As $e=1$ for $\Q_p$ and $X_0(p)$ admits a $\Q_p$-rational point (e.g. the $\infty$-cusp), and moreover $X_0(p)$ has semistable reduction over $\Z_p$ (Deligne-Rapoport \cite{DR}, V-6, or Mazur \cite{Mazur}, II-1), we know by Cor.\ \ref{loweinj} that $(\pi_1^{ab}(X_0(p)/\Q_p)^{\text{\rm geo}}_{\text{\rm ram}})^\vee$ injects to $\Phi(J_0(p))$, where $J_0(p)$ is the jacobian of $X_0(p)$. By Mazur-Rapoport \cite{MR}, we know that:

\begin{proposition}[ \cite{MR}, Theorem (A.1), b)]
$\Phi(J_0(p))$ is a cyclic group of order equal to the numerator of $\frac{p-1}{12}$.
\epr

Because the Shimura covering (\cite{Mazur}, Cor.\ (2.3)) gives a cyclic \'etale covering of $X_0(p)/\Q_p$ of order equal to the numerator of $\frac{p-1}{12}$, which is by definition completely ramified over one of the component of the special fiber, we have:

\bpr \label{X0pram}
$\pi_1^{ab}(X_0(p)/\Q_p)^{\text{\rm geo}}_{\text{\rm ram}}$ is isomorphic to $\Phi(J_0(p))$.
\epr

\brem
{\rm Similar results hold for the modular curves $X_H(p)$ between $X_1(p)$ and $X_0(p)$ corresponding to any subgroup $H\subset (\Z/p\Z)^\times/\{\pm 1\}$ at least for the prime-to-$p$ part, by the results in Conrad-Edixhoven-Stein \cite{CES}.}
\erem

\subsection{The part $\pi_1^{ab}(X_0(p)/\F_p)^{\text{\rm geo}}$}

We denote by $X_0(p)/\Z_p$ the minimal regular model of $X_0(p)/\Q_p$ following  \cite{Mazur}, and denote the special fiber by $X_0(p)/\F_p$. Here we determine the part $\pi_1^{ab}(X_0(p)/\F_p)^{\text{\rm geo}}$, thereby completing the proof of Th.\ \ref{X0pcomplete}. The result is:

\bpr \label{X0psp}
$\pi_1^{ab}(X_0(p)/\F_p)^{\text{\rm geo}}$ is a free module of over $\hat{\Z}$ with the rank equal to $\frac{g+h-1}{2}$, where $g$ is the genus of $X_0(p)$ and $h$ is the number of the supersingular points defined over $\F_p$.
\epr

\bprf
As each of the component of the special fiber $X_0(p)_{\F_p}$ of $X_0(p)$ is (geometrically) isomorphic to ${\mathbb P}^1$, we know that $\Pic^0(X_0(p)_{\F_p})$ is a torus with the character group canonically isomorphic (as the $G_{\F_p}$-module) to the first homology group $H_1(\Gamma,\Z)$ of the graph $\Gamma$ of $X_0(p)_{\ol{\F}_p}$ (\cite{SGA7} Expos\'e IX, 12.3). This graph is described in Mazur-Rapoport \cite{MR}, \S 3. In particular, if we denote the genus of $X_0(p)$ by $g$, the total number of supersingular points are $g+1$, and $H_1(\Gamma,\Z)$ is a free module over $\Z$ of rank $g$.

Now by (\ref{pi1geoasH1}), (\ref{Q/ZtoN}), and Lemma \ref{PicN}, we have:
\begin{align*}
\pi_1^{ab}(X_0(p)/\F_p)^{\text{\rm geo}} &\cong \big(\bigcup_N\Hom_{\F_p\text{\rm -gp}}(\mu_N,\Pic^0(X_0(p)_{\F_p})[N])\big)^\vee\\
&\cong \big(\bigcup_N\Hom_{G_{\F_p}}(H_1(\Gamma,\Z/N\Z),\Z/N\Z)\big)^\vee\\
&\cong \invlim{N}\big(\Hom(H_1(\Gamma,\Z/N\Z)_{G_{\F_p}},\Q/\Z)^\vee\big)\\
&\cong \invlim{N}H_1(\Gamma,\Z/N\Z)_{G_{\F_p}}\\
&\cong \invlim{N}\big(H_1(\Gamma,\Z)_{G_{\F_p}}\otimes \Z/N\Z\big)\cong H_1(\Gamma,\Z)_{G_{\F_p}}\otimes\hat{\Z}
\end{align*}
where $H_1(\Gamma,\Z)_{G_{\F_p}}$ is the $G_{\F_p}$-coinvariant of $H_1(\Gamma,\Z)$. The last isomorphism follows from the fact that $H_1(\Gamma,\Z)_{G_{\F_p}}$ is clearly a finitely generated $\Z$-module. Now it only remains to determine the group $H_1(\Gamma,\Z)_{G_{\F_p}}$ explicitly.

As we know that the vertices of $\Gamma$ is fixed by $G_{\F_p}$ as the irreducible components are defined over $\F_p$, and the edges of $\Gamma$ which correspond to the supersingular points of $X_0(p)$ are fixed or interchanged by pairs by the Frobenius automorphism of $G_{\F_p}$, according to whether the field of definition of the corresponding supersingular point is $\F_p$ or $\F_{p^2}$ (these possibilities can be visibly read off from the Table 6 of  \cite{AntwerpIV}). 

Therefore the Frobenius acts on the basis of the homology group formed by the pair of the edges corresponding to the pair of supersingular points defined over $\F_{p^2}$ by multipication by $-1$, which means that when we take the $G_{\F_p}$-coinvariant, exactly these bases vanish. Therefore if we denote the number of the supersingular points defined over $\F_p$ and the pairs of supersingular points defined over $\F_{p^2}$ respectively by $h$ and $j$ (therefore $g+1=h+2j$), the $G_{\F_p}$-coinvariant of $H_1(\Gamma,\Z)$ is a free module of rank $g-j=\frac{g+h-1}{2}$. (This rank is equal to the rank defined in the general setting in Saito \cite{SSaito}, II-Def.\ 2.5.)
\eprf

\hspace*{\parindent}

\end{document}